\documentstyle[12pt]{article}

\newcommand{\sect}[1]{\setcounter{equation}{0}
\renewcommand{\thesection}{\Roman{section}}\section{#1}}

\def\noi{\noindent}
\def\ar{\begin{array}{rcl}}
\def\an{\end{array}}
\def\ra{\rangle}

\newcommand{\eq}{\begin{equation}}
\newcommand{\eqa}{\begin{eqnarray}}
\newcommand{\en}{\end{equation}}
\newcommand{\ena}{\end{eqnarray}}
\def\1{{\bf 1}}

\def\F{\mbox{$\cal F$\,}}

\def\Fp{{\cal F}_{21}}
\def\Fu{{\cal F}^{(1)}}
\def\Fd{{\cal F}^{(2)}}

\def\Hil{{\cal H}}
\def\ot{\otimes}
\def\id{\mbox{id}}
\def\ie{\mbox{\it i.e.\/ }}
\def\eg{\mbox{\it e.g.\/ }}
\def\C{{\cal C}}
\def\up{\uparrow}
\def\idA{{\bf 1}_{\cal A}}

\def\g{\mbox{\bf g\,}}
\def\uqg{\mbox{$U_h{\/\mbox{\bf g}}$ }}

\newcommand{\tr}{\triangleright_h\,}
\newcommand{\trc}{\triangleright}
\newcommand{\tl}{\triangleleft_h\,}

\def\R{\mbox{$\cal R$\,}}

\def\A{\mbox{$\cal A$}}

\def\F{\mbox{$\cal F$}}

\def\z{\hspace*{9mm}}
\def\x{\hspace{3mm}}

\newcommand{\cn}{{\bf C}}
\newcommand{\rn}{{\bf R}}

\newcommand{\nn}{{\bf N}}
\newtheorem{prop}{Proposition}
\newtheorem{lemma}{Lemma}
\newtheorem{theorem}{Theorem}

\begin{document}

\begin{titlepage}
\begin{center}
~
~

\vskip.6in

{\Large \bf Deforming Maps for Lie Group Covariant
             Creation \& Annihilation Operators}

\vskip.4in

Gaetano Fiore* 

\vskip.25in

{\em Sektion Physik der Ludwig-Maximilians-Universit\"at
M\"unchen\\
Theoretische Physik --- Lehrstuhl Professor Wess\\
Theresienstra\ss e 37, 80333 M\"unchen\\
Federal Republic of Germany}
\end{center}
\vskip1in
\begin{abstract}
Any deformation of a Weyl or Clifford algebra ${\cal A}$ can be
realized through a `deforming map', i.e. a formal change of
generators in ${\cal A}$. This is true in 
particular if ${\cal A}$ is covariant
under a Lie algebra $\g$  and its deformation is induced by
some triangular deformation $U_h\g$ of the Hopf algebra 
$U\g$. We propose a systematic method to 
construct all the corresponding deforming
maps, together with the corresponding realizations
of the action of $U_h\g$. The method is then generalized
and explicitly applied to the case that $U_h\g$ is the quantum
group ${\cal U}_hsl(2)$. A preliminary study of the status of
deforming maps at the representation level 
shows in particular that 
`deformed' Fock representations induced by a compact $U_h\g$ 
can be interpreted as standard `undeformed' Fock representations 
describing particles with ordinary Bose
or Fermi statistics.
\end{abstract}
\vfill
\begin{center}
\end{center}
\noi \hrule
\vskip.2cm
\noi{\footnotesize *EU-fellow, TMR grant ERBFMBICT960921.
\qquad {\it e-mail: }Gaetano.Fiore \ @ \ physik.uni-muenchen.de}
\end{titlepage}
\newpage
\setcounter{page}{1}

\sect{Introduction}

In recent years the idea of noncocommutative Hopf
algebras \cite{abe}
(in particular quantum groups \cite{dr2})
as candidates for generalized symmetry
transformations in quantum physics 
has raised an increasing interest. One way to
implement this idea in quantum field theory
or condensed matter physics would be to deform
the canonical commutation relations (CCR) of
some system of mode creators/annihilators, covariant
under the action of a Lie algebra $\g$,
in such a way that they become covariant under the action
of a noncocommutative deformation $U_h\g$
(with deformation parameter $h$)
of the cocommutative Hopf algebra $U\g$, as it has been
done \eg in Ref. \cite{puwo,wezu} for the 
${\cal U}_hsl(N)$ covariant Weyl algebra in $N$ dimensions. 

As a toy model for
these deformations one can consider
the deformed Weyl algebra $\A_h$
in 1 dimension \cite{mac}
with generators fulfilling the `quantum' commutation relation
(QCR)
\eq
\tilde A\,\tilde A^+=\1+ q^2\tilde A^+\,\tilde A
\label{proto}
\en
with $q=e^h$. 
When $q=1$ the above
reduces to the classical Weyl algebra $\A$
\eq
a \,a^+=\1+a^+a .
\en
If we define $n :=a^+a $, $(x)_z={z^x-1\over z-1}$
and \cite{zachos3}
\eq
A:=a \sqrt{\frac{(n )_{q^2}}{n }}
\qquad\qquad
A^+:=\sqrt{\frac{(n )_{q^2}}{n }}a^+,
\label{traproto}
\en
we find out that $A,A^+$ fulfil the QCR
(\ref{proto});  hence we can define an 
algebra homomorphism $f:\A_h\rightarrow \A[[h]]$,
or ``deforming map'' (in the terminology of
Ref. \cite{zachos1,zachos2}), starting from
\eq
f(\tilde A)=A\qquad\qquad
f(\tilde A)=A^+.
\en
The RHS(\ref{traproto}) have to be
understood as formal power series in the
deformation parameter $h$.

We are interested in
deformed multidimensional Weyl or Clifford algebras 
$A_h$ where the QCR:

\begin{enumerate}

\item keep a quadratic
structure as in eq. (\ref{proto}), so that one can represent
the generators as creation or annihilation operators;

\item are covariant under the action $\tilde\tr$ of $U_h\g$.

\end{enumerate}

More precisely, the generators $\tilde A^+_i,\tilde A^i$ should 
transform linearly under the action of $\tilde\tr$,
\eqa
x\,\tilde\tr\tilde A^+_i &=& \tilde{\rho}^j_i(x)A^+_j\\
x\,\tilde\tr\tilde A^i &=& \tilde{\rho}^\vee{}^i_j(x)A^j
\label{melma}
\ena
with $\tilde{\rho}^\vee$ being the the contragradient
representation of the representation $\tilde\rho$ of
$U_h\g$ ($\tilde{\rho}^\vee=\tilde{\rho}^T\circ S_h$, where
$S_h$ is the antipode of $U_h\g$ and $^T$ is the operation of
matrix transposition).

In this work we essentially stick to the case that $U_h\g$
is {\it triangular}; we treat the general quasitriangular case
in Ref. \cite{fiocomp}. In the former case one can show easily
that, for arbitrary $\tilde{\rho}$, $U_h\g$-covariant QCR
are given by
\eqa
\tilde A^i\tilde A^j          & = & \pm R^{ij}_{vu}\tilde A^u\tilde A^v \\
\tilde A^+_i\tilde A^+_j & = & \pm R_{ij}^{vu}\tilde A^+_u\tilde A^+_v
\label{QCR0}\\
\tilde A^i\tilde A^+_j     & = & \delta^i_j\idA \pm R^{ui}_{jv}
\tilde A^+_u\tilde A^v;
\ena
here the sign $\pm$
refers to the Weyl/Clifford case respectively, and $R$ is the
corresponding `$R$-matrix' of $U_h\g$\footnote{One just
has to note that
$R=\tilde{\rho}\ot\tilde{\rho}\R$, where $\R$ is the universal
triangular structure of $H_h$, and that 
$\tau\circ\Delta_h(x)=\R\Delta_h(x)\R^{-1}$ ($\tau$ denotes the
flip operator).}.

$A_h$ is a left-module algebra of $U_h\g$:
 the `quantum' action $\tilde\tr$
is extended to products of the generators as a left-module 
algebra map $\tilde\tr:U_h\times \A_h\rightarrow \A_h$
(\ie consistently with the QCR) using the coproduct
$\Delta_h(x)=\sum_{\mu} x^{\mu}_{(\bar 1)}\ot x^{\mu}_{(\bar 2)}$
of $U_h\g$,
\eq
x\,\tilde\tr(a\cdot b)=
\sum_{\mu} (x^{\mu}_{(\bar 1)}\,\tilde\tr a)\cdot 
(x^{\mu}_{(\bar 2)}\,\tilde\tr b),
\label{gigi}
\en
because the
(\ref{QCR0}) are covariant under (\ie compatible with) $\tilde\tr$.

The existence of deforming maps for arbitrary 
(\ie not necessarily of the kind described above)
deformations of Weyl (or Clifford) algebras
is a consequence \cite{gerst}
of a theorem \cite{ducloux}
asserting the triviality
of the cohomology groups of the latter
(see Ref. \cite{mathias,fiocomp} for
an effective and concise presentations of these results.
See also Ref. \cite{flato}, where the problem
of stability of quantum mechanics under deformations
was addressed for the first time.).
However, no general method for their
explicit construction is available. 
Actually, using cohomological arguments,
one can also easily show that deforming
maps are unique up to a inner automorphism,
\eq
f\rightarrow f_{\alpha}:=\alpha f(\cdot) \alpha^{-1}
\qquad\qquad \alpha=\idA+ O(h);
\label{auto}
\en
therefore it is enough to construct one to
find all of them.

In this work we present a general method which allows,
given a triangular Hopf algebra $U_h\g$ and any 
$U_h\g$-covariant deformed Weyl or Clifford algebra,
to explicitly construct the corresponding deforming maps $f$
and the corresponding realizations $\tr$ of $\tilde\tr$
[ $\tr$ is defined by $\tr:=(\id\ot f)\circ\tilde\tr\circ(\id\ot f^{-1})$].
In a first attempt to generalize our construction procedure
to {\it quasitriangular} $U_h\g$, we also generalize the constrution to the
case that $U_h\g$ is the quantum group ${\cal U}_hsl(2)$
and $\tilde \rho$ is its fundamental representation.
Finally we investigate on the status of deforming maps
at the representation-theoretic level.

The construction method 
is based (Sect. \ref{preli1}) on use of the 
Drinfel'd-Reshetikhin twist $\F$ \cite{dr1,resh}, intertwining
between the coproducts of $U\g$ and $U_h\g$, and on the fact that 
within ${\cal A}[[h]]$ one can realize both the action $\trc$
of $U_h\g$ (Section  \ref{preli1}) and the action $\tr$ in an
`adjoint-like' way. We show first (Section
\ref{qcov}) that $\F$ can be used in a {\it universal}
way to construct, within ${\cal A}_h$, 
$U_h\g$-tensors out of $U\g$-tensors, and in particular
out of $a^+_i,a^i$ objects $A^+_i,A^i$ that transform
under $\tr$ as in formula (\ref{melma}).
Then (Section \ref{trcr}) we verify that the objects
$A^+_i,A^j$ really satisfy the QCR (\ref{QCR0}). In
Section \ref{sl2cr} we generalize our construction 
(by means of the Drinfel'd twist \cite{dr3}) to the
case of deformed Weyl \& Clifford algebras with generators
belonging to the fundamental
representation of the quantum group ${\cal U}_hsl(2)$;
the deforming map is again completely explicit
thanks to the semiuniversal
expression \cite{zachos2} for $\F$. We compare our
deforming map with the one previously found in 
Ref. \cite{oleg}. At the representation-theoretic
level it would be natural to interpret deforming
maps as ``operator maps'', in other words as 
intertwiners between the representations of
$\A$ and $\A_h$. However we have to expect that,
in the role of intertwiners, deforming maps
may become singular at $h=0$, because the representation
theories of $\A$, $\A_h$ are in general rather different.
In Section \ref{repth} we show that there is always a 
$*$-representation of $\A_h$ which is intertwined
by $f$ with the Fock representation of $\A$;
this allows to interpret $\tilde A^i,\tilde A^+_i$ as
`composite' operators on a classical Fock space
describing ordinary Bosons and Fermions. We
also explicitly show that $f_{\alpha}^{-1}$ is
ill-defined as an intertwiner from the remaining
(if any) unitarily inequivalent $*$-representations
of $\A_h$.

On the basis of the above result we conclude that also at the
quantization-of-field level noncommutative Hopf algebra
symmetries are not necessarily incompatible with Bose or
Fermi statistics (contrary to what is often claimed).
We arrived at the same conclusion at the first-quantization
level in Ref. \cite{fioschu,fiopro}, where the initial
motiation for the present work has originated.
The connection between the two approaches through second
quantization will be described elsewhere.

\sect{Preliminaries and notation}

\subsection{Twisting groups into quantum groups}
\label{preli1}

Let $H =(U\g,m,\Delta , \varepsilon,S )$ be the
cocommutative Hopf algebra associated to the universal
enveloping (UE) algebra $U\/ \g$ of a Lie algebra \g.
The symbol $m$ denotes the multiplication (in the
sequel it will be dropped in the obvious way
$m(a\ot b)\equiv ab$, unless explicitly required), whereas
$\Delta , \varepsilon,S $  the comultiplication, counit
and antipode respectively.

Let $\F\in U\g[[h]]\ot U\g[[h]]$ (we will write
$\F=\F^{(1)}\ot \F^{(2)}$, in a Sweedler's
notation with {\it upper} indices; in the RHS a sum
$\sum_i \F^{(1)}_i\ot \F^{(2)}_i$
of many terms is implicitly understood) be
a `twist', \ie an element satisfying the relations
\eqa
&&(\varepsilon\ot \id)\F=\1=(\id\ot \varepsilon)\F
\label{cond2}\\
&&\F=\1\ot\1 +O(h)
\label{cond2bis}
\ena
($h\in\cn$
is the `deformation parameter', and $\1$ the unit in $U\g$; from the second
condition it follows that $\F$ is invertible as a power series).
It is well known \cite{dr1} that if \F \/ also satisfies the relation
\eq
(\F\ot \1)[(\Delta \ot \id)(\F)]=(\1 \ot \F)[(\id \ot \Delta )(\F),
\label{cond1}
\en
and $(U_h\g,m_h)$ is an algebra
isomorphic to  $U\g[[h]]$ with isomorphism, say,
$\varphi_h:U_h\g\rightarrow U\g[[h]]$ [in particular, if $U_h\g=U\g[[h]]$]
and $\varphi_h=\id\x (\mbox{mod } h)$, or even $\varphi_h=\id$],
then one can construct a
triangular non-cocommutative Hopf algebra
$H_h=(U_h\g,m_h,\Delta_h,\varepsilon_h,S_h,\R)$
having an isomorphic  algebra
structure  [$m_h= \varphi_h^{-1}\circ m\circ (\varphi_h\ot\varphi_h)$],
an isomorphic counit
$\varepsilon_h:=\varepsilon\circ \varphi_h^{-1}$,
comultiplication and antipode defined by
\eq
\Delta_h(a)=(\varphi_h^{-1}\ot\varphi_h^{-1})\big
\{\F\Delta [\varphi_h(a)]\F^{-1}\big\},\qquad\qquad\qquad
S_h(a)=\varphi_h^{-1}\{\gamma^{-1} S [\varphi_h(a)]\gamma\},
\label{def1}
\en
where
\eq
\gamma := S \F^{-1(1)}\cdot\F^{-1(2)}, \z\x \qquad\qquad
\gamma^{-1} = \F^{(1)}\cdot S \F^{(2)},
\label{def2}
\en
and (triangular) universal R-matrix
\eq
\R:=[\varphi_h^{-1}\ot\varphi_h^{-1}](\Fp\F^{-1}), \qquad\qquad\qquad
\F_{21}:=\F^{(2)}\ot\F^{(1)}.
\label{rmat}
\en
Condition (\ref{cond1}) ensures that $\Delta_h$ is
coassociative as $\Delta $.
The inverse of $S_h$ is given by
$S_h^{-1}(a)=\varphi_h^{-1}\{\gamma'S [\varphi_h(a)]\gamma'^{-1}\}$,
where
\eq
 \gamma' := \F^{(2)}\cdot S \F^{(1)} \z\x \qquad\qquad
\gamma'^{-1} = S \F^{-1(2)}\cdot\F^{-1(1)};
\label{def2'}
\en
$\gamma^{-1}\gamma'\in\mbox{Centre(U\g)}$, and
$S \gamma=\gamma'^{-1}$.

Conversely, given a
$h$-deformation $H_h=(U_{h}\g,m_h,\Delta_h,\-\varepsilon_h,S_h,\R)$
of $H $ in the form of a
triangular Hopf algebra, one can find \cite{dr1} and an isomorphism
$\varphi_h:U_{h}\rightarrow U\g[[h]]$
an invertible \F ~satisfying conditions
(\ref{cond2}), (\ref{cond2bis}), (\ref{cond1}) such that $H_h$ can be obtained
from $H $  through formulae
(\ref{def1}),(\ref{def2}),(\ref{def2'}).

Examples of \F's satisfying conditions  (\ref{cond1}),
(\ref{cond2}), (\ref{cond2bis}) are provided \eg by the
socalled `Reshetikhin twists' \cite{resh}
\eq
\F:= e^{h \omega_{ij}h_i\ot h_j},
\label{reshet}
\en
where $\{h_i\}$ is a basis of the Cartan subalgebra of \g
and $\omega_{ij}=-\omega_{ji}\in\cn$. A less obvious example
is for instance the `Jordanian' deformation of
Ref. \cite{ohn}.
   
A similar result to the above holds for genuine quantum groups.
A well-known theorem by Drinfel'd, Proposition 3.16
in Ref. \cite{dr3} proves, for any quasitriangular
deformation $H_h=(\uqg,m_h,\Delta_h,\varepsilon_h,S_h,\R)$
\cite{dr2,frt}
of $U\g$, with  \g a simple finite-dimensional Lie algebra,
the existence of an algebra  isomorphism
$\varphi_h: \uqg \rightarrow U\g[[h]]$ and an
invertible \F ~satisfying condition
(\ref{cond2}) such that $H_h$ can be obtained from $H $
through formulae
(\ref{def1}),(\ref{def2}),(\ref{def2'}),  as well, after
identifying $h=\ln q$.
This \F~ does not satisfy
condition (\ref{cond3}), however the (nontrivial) coassociator
\eq
\phi:=[(\Delta \ot \id)(\F^{-1})](\F^{-1}\ot \1)(\1 \ot \F)
[(\id \ot \Delta )(\F)
\label{defphi}
\en
still commutes with $\Delta^{(2)} (U\g)$,
\eq
[\phi,\Delta^{(2)} (U\g)]=0,
\label{ginv}
\en
 thus explaining why
$\Delta_h$ is coassociative in this case, too. The corresponding
universal (quasitriangular) R-matrix \R is related to
\F ~by
\eq
\R=[\varphi_h^{-1}\ot\varphi_h^{-1}](\Fp q^{t\over 2}\F^{-1}),
\label{defR}
\en
where $t:=\Delta (\C)-\1\ot \C-\C\ot \1$ is
the canonical invariant element
in $U\g\ot U\g$ ($\C$
is the quadratic Casimir). The twist $\F$ is
defined (and unique) up to the transformation
\eq
\F\rightarrow \F T,
\label{trasfo0}
\en
where $T$ is a $\g$-invariant
[\ie commuting with $\Delta (U\g)$] element
of $U\g[[h]]^{\ot^2}$ such that
\eq
T=\1\ot\1+O(h),\z\z
(\varepsilon\ot\id) T=\1=(\id\ot\varepsilon) T.
\en
A function
\eq
T=T\left(\1\ot {\cal C}_i, {\cal C}_i \ot \1, \Delta ({\cal C}_i)\right)
\label{hil14}
\en
of the Casimirs ${\cal C}_i\in U\g$
of $U\g$ and of their coproducts clearly is \g-invariant.

In general, as a consequence of the existence of an isomorphism $\varphi_h$,
representations $\tilde\rho,\rho$ of deformed and undeformed algebrae
are in one-to-one correspondence (except for special values of $h$ making
it singular) through
\eq
\rho=\tilde\rho\circ\varphi_h.
\en

\medskip
A special case of interest is when
$U\g$ is a $*$-Hopf algebra and \F ~ is unitary, 
\eq
\F^{*\ot *}=\F^{-1};
\label{cond4}
\en
note that in this case
\eq
\gamma'=\gamma^*.
\label{stuff}
\en
One can show \cite{jurco} that \F~
can always be made unitary if \g\/ is compact.

We will often use a `tensor notation' for our formulae:
eq. (\ref{cond1}) will read
\eq
\F_{12}\F_{12,3}=\F_{23}\F_{1,23},
\label{cond3}
\en
and definition (\ref{defphi})
$\phi\equiv\phi_{123}=\F_{12,3}^{-1}\F_{12}^{-1}\F_{23}\F_{1,23}$,
for instance; the comma separates the tensor factors 
{\it not} stemming from the coproduct.
For practical purposes it will be often convenient in the sequel to use
 the Sweedler's notation with {\it lower} indices
$\Delta (x)\equiv x_{(1)}\ot x_{(2)}$ for the
cocommutative coproduct
(in the RHS a sum
$\sum_i x^i_{(1)}\ot x^i_{(2)}$
of many terms is implicitly understood); similarly, we will use the
Sweedler's notation
$\Delta ^{(n-1)}(x) \equiv x_{(1)}\ot\ldots\ot x_{(n)}$
for the $(n\!-\!1)$-fold coproduct.
For the non-cocommutative coproducts $\Delta_h$, instead, we will
use a Sweedler's notation with barred indices:
$\Delta_h(x)\equiv x_{(\bar 1)}\ot x_{(\bar 2)}$.

\subsection{Classical $U\g$-covariant creators and annihilators}
\label{preli2}

Let \A~ be the unital algebra  generated by $\idA$ and
elements $\{a^+_i\}_{i\in I}$ and $\{a^j\}_{j\in I}$
satisfying the (anti)commutation relations
\eqa
[a^i\, , \,a^j]_{\pm}            &=& 0 \nonumber\\~
[a^+_i\, , \, a^+_j]_{\pm}  &=& 0 \label{ccr}\\~
[a^i\, , \, a^+_j]_{\pm}       &=& \delta_j^i\idA
\nonumber
\ena
(the $\pm$ sign denotes commutators and
anticommutators respectively),
belonging respectively to some
representation $\rho$ and to its
contragradient $\rho^\vee = \rho^T \circ S $ of $H $
(${}^T$ is the transpose):
\eq
\begin{array}{rcl}
x\trc a^+_i &=&\rho(x)^l_ia^+_l  \cr
x\trc a^i&=&\rho(S x)^i_la^l \cr
\end{array}
\qquad\qquad\qquad x\in U\g,\z \rho(x)^i_j\in \cn.
\label{covl}
\en
Equivalently, one says that $a^+_i,a^i$ are ``covariant''
under $\trc$, or that they span two (left) modules of $U\g$:
\eq
(xy)\trc a= x\trc(y\trc a),\z\z\z x,y\in U\g,
\label{chicca}
\en
with either $a=a^i$ or $a^+_i$.

\A~ is a (left) module algebra of $(H ,\trc)$, if the action
$\trc$ is extended on the
whole \A~ by means of the (cocommutative) coproduct:
\eq
x\trc (ab)=(x_{(1)}\trc a)(x_{(2)}\trc b).
\en
Then property (\ref{chicca}) holds for all $a\in\A$.

Setting
\eq
 \sigma(X):=
\rho(X)^i_ja^+_i a^j
\label{jordan}
\en
for all $X\in \g$, one  finds that
$\sigma: \g\rightarrow \A$
is a Lie algebra homomorphism,
so that  $\sigma$ can be extended to
all of $U\g$ as an algebra
homomorphism $\sigma: U\g\rightarrow \A$;
on the unit element we set
$\sigma(\1_{U\g}):=\idA$. $\sigma$ can be seen
as the generalization of the Jordan-Schwinger
realization of $\g=su(2)$ \cite{bied} [formula (\ref{homo})].

Then it is easy to check the following
        
\begin{prop}
The (left) action $\trc:U\g\times \A\rightarrow \A$ can be
realized in an `adjoint-like' way:
\eq
x\trc a=\sigma(x_{(1)})\: a\: \sigma(S  x_{(2)}) ,
\qquad\qquad\qquad x\in U\g, \z a\in\A.
\label{cov}
\en
\label{lemma}
\end{prop}

In the specially intersting case of a compact section $\g$ 
(with $*$-structure ``$*$'') one can
introduce 
in \A~ a $*$-structure, the `hermitean conjugation'
(which we will denote by ${}^{\star}$),  
such that
\eq
(a ^i)^{\star}= a^+_i.
\label{star1}
\en
Correspondingly, $\rho$ is a $*$-representation
($\rho\circ \star=*\circ\rho^T$) and
$\sigma$ becomes a $*$-homomorphism, \ie
$\sigma\circ *=\star\circ \sigma$.

\sect{Quantum covariant creators and annihilators}
\label{qcov}

Let $H_h$ and $\varphi_h$ be as in section \ref{preli1}.
Clearly, $\sigma_{\varphi_h}:=\sigma\circ\varphi_h$
is an algebra homomorphism 
$\sigma_{\varphi_h}:U_h\g\rightarrow \A[[h]]$.
Inspired by proposition \ref{prop1} we are led
to define
\eq
x\tr a := \sigma_{\varphi_h}(x_{(\bar 1)}) a 
\sigma_{\varphi_h}(S_h x_{(\bar 2)}).
\label{gog}
\en
Using the Hopf algebra
axioms it is straightforward to prove the relations
[cfr. relations (\ref{gigi})]
\eq
\begin{array}{rcl}
(xy) \tr a & = & x \tr (y\tr a) \cr
x \tr  (ab)  & = &(x_{(\bar 1)}\tr a)
(x_{(\bar 2)}\tr b),
\end{array}
\z\z
\forall x,y\in U\g[[h]],\x \forall a,b\in\A[[h]].
\en
In other words
\begin{prop}
The definition (\ref{gog}) realizes $\tilde\tr$
(the left action of $H_h$) on the algebra
$\A[[h]]$.
\end{prop}

However, $a^+_i,a ^j$ are {\it not} covariant 
w.r.t. to $\tr$. One may ask whether there exist some
objects $A^+_i,A^j\in \A$ that {\it are} covariant
under $\tr$ and transform as in eq. (\ref{QCR0}).

The answer comes from the crucial
\begin{prop}
The elements
\eqa
A_i^+ &:= &\sigma(\Fu)a_i^+
\sigma(S \Fd\gamma)\in\A[[h]] \label{def3.1}\\
A^i&:= &\sigma(\gamma'S \F^{-1(2)})A^i 
\sigma(\F^{-1(1)})\in \A[[h]]                          
\label{def3.2}
\ena
are ``covariant'' under $\tr$, more precisely 
belong respectively to the representation
$\tilde\rho$ and to its quantum contragredient
$\tilde\rho^\vee=\tilde\rho^T\circ S_h$ of $U_h\g$ acting 
through $\tr$:
\eq
x\tr A^+_i=\tilde\rho(x)^l_iA^+_l\z\z\z x\tr A^i=\tilde\rho(S_hx)^i_mA^m.
\en 
\label{prop1}
\end{prop}
{\it Proof.} Due to relation (\ref{def1}), $\F$ is an intertwiner 
between $\Delta_h$ and $\Delta $ (in this proof 
we drop the symbol $\varphi_h$):
\eq
x_{(\bar 1')}\Fu\ot x_{(\bar 2')}\Fd=\Fu x_{(1')}\ot \Fd x_{(2')}.
\label{inter}
\en
Applying $\id\ot S $ on both sides of the equation and 
multiplying the result by $\1\ot\gamma$ from the right
we find [with the help of relation (\ref{def2})]
\eq
x_{(\bar 1')}\Fu\ot (S \Fd) \gamma S_h x_{(\bar 2')}=\Fu x_{(1')}\ot
(S  x_{(2')})(S \Fd)\gamma.
\en
Applying $\sigma\ot\sigma$ to both sides and sandwiching
$a^+_i$ between the two tensor factors we find
\eq
\sigma(x_{(\bar 1')})A^+_i \sigma(S_h x_{(\bar 2')})=
\sigma(\Fu)\sigma( x_{(1')})
a^+_i\sigma(S x_{(2')})\sigma[(S \Fd)\gamma],
\en
which, in view of formula (\ref{gog}),
proves the first relation. 

To prove the second relation, let us note that 
relation (\ref{def1})  implies an analogous
relation
\eq
\Delta_h(a)\widetilde{\F}=\widetilde{\F}\Delta (a),
\qquad\qquad\qquad
\mbox{with~}\widetilde{\F}:= [\gamma' S  \F^{-1~(2)}   
\ot  \gamma'  S  \F^{-1~(1)} ]\Delta (S \gamma) .
\label{def1'}
\en
This can be shown by applying in the order 
the following operations to both sides of
eq. (\ref{def1}): multiplying by $\F^{-1}$ from 
the left and from the right, applying $S \ot S $,
multiplying by $\gamma'\ot \gamma'$
from the left and by $\Delta (S \gamma)$
from the right, replacing $a\rightarrow S_h x$,
using the properties (\ref{def1}) and
$(S_h\ot S_h)\circ \Delta_h=\tau \circ \Delta_h \circ S_h$.
Next, we observe that $A^i$ can be rewritten as
\eq
A^i=\sigma[\widetilde{\F}^{(1)}S (\gamma^{-1})_{(1)}]a^i
\sigma[(\gamma^{-1})_{(2)}S \widetilde{\F}^{(2)}\gamma]=
\sigma(\widetilde{\F}^{(1)})a^l
\sigma(S \widetilde{\F}^{(2)}\gamma)\rho(\gamma^{-1})^i_l;
\label{pippo}
\en
whence, reasoning as for the first relation,
\eqa
\sigma(x_{(\bar 1')})A^i \sigma(S_h x_{(\bar 2')}) 
&\stackrel{(\ref{pippo})}{=} & \sigma(\widetilde{\F}^{(1)})\sigma( x_{(1')})
a^l \sigma(S x_{(2')})\sigma[(S \widetilde{\F}^{(2)})\gamma]
\rho(\gamma^{-1})^i_l                                   \nonumber    \cr
~ & \stackrel{(\ref{cov})}{=}&
\sigma(\widetilde{\F}^{(1)})a^l \sigma[(S \widetilde{\F}^{(2)})\gamma]
\rho(\gamma^{-1}S  x)^i_l                           \nonumber   \cr
~ & \stackrel{(\ref{def1})}{=}&\sigma(\widetilde{\F}^{(1)})a^l
\sigma[(S \widetilde{\F}^{(2)})\gamma]
\rho(S_hx \cdot\gamma^{-1})^i_l                      \nonumber    \cr
~ &\stackrel{(\ref{pippo})}{=} & \rho(S_hx)^i_lA^l
\ena
which proves the second relation.  $\Box$

\medskip

{\it Remark 1} The proposition clearly
holds even if one chooses in 
formulae (\ref{def3.1}), (\ref{def3.2}) two \F's differing
by a tranformation \ref{trasfo0}. We shall exploit this freedom
when $U_h\g$ is genuinely quasitriangular.

{\it Remark 2} Note that in the $*$-Hopf algebra case 
eq. (\ref{star1}), (\ref{cond4}), (\ref{stuff}) imply
\eq
(A^i)^{\dagger}= A^+_i.
\label{star2}
\en

{\it Remark 3} Under the right action $\tl$
($a\tl x:= (S_h^{-1}x)\tr a$ with $a\in\A,\x x\in U\g$)
the covariance properties of $A^i,A^+_i$ read
\eq
A^i\tl x
=\rho(x)^i_lA^l\z\z\z  A^+_i\tl x=\rho(S_h^{-1}x)_i^mA^+_m.
\en 

{\it Remark 4} $\sigma_{\varphi_h}$ is not the only algebra
homomorphism $U_h\g\rightarrow \A[[h]]$. For any
$\alpha=\idA+O(h)$ we find a new one by setting
\eq
\sigma_{\varphi_h,\alpha}(x)=\alpha \sigma_{\varphi_h}(x)\alpha^{-1};
\en
correspondingly, we can define a new realization of $\tilde\tr$
by
\eq
x\trc_{h,\alpha} a := \sigma_{\varphi_h,\alpha}(x_{(\bar 1)}) a 
\sigma_{\varphi_h,\alpha}(S_h x_{(\bar 2)}).
\label{goga}
\en
Covariant objects under $\tr_{\alpha}$ will be given by
\eq
A^+_{i,\alpha}:=\alpha A^+_i\alpha^{-1}, \z\z\z
A^{i,\alpha}:=\alpha A^i\alpha^{-1}.
\label{fuffi}
\en
If relations (\ref{star2}) hold, we can preserve them by
choosing $\alpha^{\star}=\alpha^{-1}$.

To conclude this section, let us give useful alternative 
expressions for $A^+_i,A^i$ by `moving' to  the right/left 
  past $a^+_i,a^i$ the expressions
$\sigma (\cdot)$ lying at their  left/right in definitions 
(\ref{def3.1}), (\ref{def3.2}). 

\begin{lemma}
If ${\cal T}\in U\g[[h]]^{\ot^3}$
is \g-invariant (\ie $[{\cal T},U\g[[h]]^{\ot^3}]=0$) then
$m_{ij}S_{i}{\cal T}$, $m_{ij}S_{j}{\cal T}$ ($i,j=1,2,3$, $i\neq j$)
are $\g$-invariants belonging to $U\g[[h]]^{\ot^2}$.
\label{lemmino}
\end{lemma}
 (Here
$S_{i}$ denotes $S $ acting on
the $i$-th tensor factor, and
$m_{ij}$ multiplication of the $i$-th tensor factor by
the $j$-th from the right.)

$Proof$. For instance,
\[
{\cal T}^{(1)}x_{(1)}\ot {\cal T}^{(2)}x_{(2)}\ot {\cal T}^{(3)}x_{(3)}\ot
x_{(4)} = x_{(1)}{\cal T}^{(1)}\ot x_{(2)} {\cal T}^{(2)}\ot
x_{(3)}{\cal T}^{(3)}\ot x_{(4)}
\z{\Rightarrow}
\]
\[
\stackrel{(m_{23})^2\circ S_{3}}{\Longrightarrow}\x\x
{\cal T}^{(1)}x_{(1)}\ot {\cal T}^{(2)}x_{(2)} S x_{(3)}S {\cal T}^{(3)}
x_{(4)} = x_{(1)}{\cal T}^{(1)}\ot x_{(2)} {\cal T}^{(2)}S {\cal T}^{(3)}
S x_{(3)} x_{(4)}
\]
for any $x\in U\g[[h]]$, whence (because of
$a_{(1)}\,S a_{(2)}=\varepsilon(a)=S a_{(1)}a_{(2)}$)
\eq
{\cal T}^{(1)}x_{(1)}\ot {\cal T}^{(2)}S {\cal T}^{(3)}x_{(2)}=
x_{(1)}{\cal T}^{(1)}\ot x_{(2)} {\cal T}^{(2)}S {\cal T}^{(3)},
\en
so that
${\cal T}^{(1)}\ot {\cal T}^{(2)}S {\cal T}^{(3)}\in U\g[[h]]^{\ot^2}$ is
$\g$-invariant. $\Box$

\begin{prop}
If $H_h$ is triangular the definitions (\ref{def3.1}),
(\ref{def3.2}) amount to
\eq
\begin{array}{rcl}
A_i^+ &= &a_l^+\sigma(\F^{-1(2)})\rho(\F^{-1(1)})_i^l\\
A_i^+ &= &\rho(S \Fu\gamma'{}^{-1})_i^l\sigma(\Fd)a^+_l\\
A^i &= &\rho(\Fu)^i_l\sigma(\Fd) a^l \\
A^i &= &a^l\sigma(\F^{-1(2)})\rho(\gamma^{-1} S \F^{-1(1)})^i_l.
\label{def4}
\end{array}
\en
If $H_h$ is quasitriangular the same formulae hold with in
general four different $\F$'s [related to each other by
transformations (\ref{trasfo0})].
\label{prop2}
\end{prop}

$Proof.$ Observing that
\eqa
\sigma(x)a&=&\sigma(x_{(1)} ) a\sigma(S x_{(2)}\cdot x_{(3)} )\\
a\sigma(x)&=&\sigma(x_{(3)} S  x_{(2)} )a\sigma(x_{(1)} )
\label{com}
\ena
for all $x\in U\g$, $a\in\A$, we find from relations
(\ref{def3.1}), (\ref{def3.2}) and (\ref{covl})
\eqa
A_i^+ &= &a_l^+\sigma(\F^{(1)}_{(2)} S \F^{(2)}\gamma)
\rho(\F^{(1)}_{(1)})_i^l\label{def5.1}\\
A^i&= &\sigma(\gamma'S \F^{-1(2)}\cdot
\F^{-1(1)}_{(2)})\rho(\F^{-1(1)}_{(1)} )^i_l A^l .
\label{def5.2}
\ena
On the other hand, 
applying the previous lemma to ${\cal T}=\phi$
[formula {\ref{defphi})] we find that
\[
T:=m_{23}(\id\ot\id\ot S )\phi\stackrel{(\ref{def2})}{=}
(\F^{-1(1)}{}_{(1)}\ot \F^{-1(1)}{}_{(2)})
\F^{-1}(\1\ot \gamma^{-1}S \F^{-1(2)})
\]
is $\g$-invariant, whence one easily finds the relation
\eq
\Fu{}_{(1)} \ot \Fu{}_{(2)}S (\Fd)\gamma
=T^{-1}\F^{-1}\equiv \F^{'-1},
\en
after noting that
$[T,\F^{\pm 1(1)}{}_{(1)}\ot \F^{\pm 1(1)}{}_{(2)}]=0$.
Replacing in eq. (\ref{def5.1}) one finds
relation (\ref{def4})$_1$. In the triangular case
$\phi\equiv \1^{\ot^3}$, implying $T\equiv \1^{\ot^2}$
anf $\F'=\F$. Similarly one proves the other relations.
 Relations
(\ref{def4})$_3$, (\ref{def4})$_4$,  can be found also
more directly starting from relations
(\ref{def4})$_1$, (\ref{def4})$_2$
by observing 
that in the unitary-\F ~case they follow from the
latter two by applying the $*$-conjugation.
$\Box$

\sect{Quantum commutation relations: the triangular case}
\label{trcr}

\begin{theorem}
If the noncocommutative Hopf algebra $H_h$
is triangular [\ie the
twist \F satisfies equation (\ref{cond2})], then
$A^i,A_j^+$ close the quadratic commutation relations
\eqa
A^iA^+_j     & = & \delta^i_j\idA \pm R^{ui}_{jv}A^+_uA^v \,
\label{qccr1}\\
A^iA^j          & = & \pm R^{ij}_{vu}A^uA^v 
\label{qccr2}\\
A^+_iA^+_j & = & \pm R_{ij}^{vu}A^+_uA^+_v
\label{qccr3} 
\ena
where $R$ is the (numerical) quantum R-matrix of
$U\g$ in the representation $\rho$,
\eq
R^{ij}_{hk}:=\big[(\rho\ot \rho)(\R)\big]^{ij}_{hk}.
\en
\end{theorem}

$Proof$. Beside eq.'s (\ref{def4}), we will need
their `inverse' relations:
\eqa
a_l^+ &= &A_i^+\sigma(\Fd)\rho(\Fu)^i_l\nonumber\\
a_l^+ &= &\rho(\gamma S \F^{-1(1)})^i_l
\sigma(\F^{-1(2)})A^+_i \nonumber\\
a^l&= &\rho(\F^{-1(1)})^l_i\sigma(\F^{-1(2)}) A^i
\label{def4in} \\
a^l&= & A^i\sigma(\Fd)\rho(S \Fu\cdot \gamma)^l_i.
\nonumber
\ena
Using eq.'s (\ref{def4}),
\eqa
A^iA_j^+ \! &\!\stackrel{(\ref{def4})}{=}\!& \! 
\rho(\Fu)^i_l\sigma(\Fd) a^l 
a_m^+\sigma(\F^{-1(2')})\rho(\F^{-1(1')})_j^m
\nonumber \\
 \! &\!\stackrel{(\ref{ccr})}{=} \! &  \! 
\rho(\Fu)^i_l\sigma(\Fd) [\delta^l_m \pm a_m^+
a^l] \sigma(\F^{-1(2')})\rho(\F^{-1(1')})_j^m
\nonumber \\
 \! &\!\!\stackrel{(\ref{com})\!,\!(\ref{covl})}{=}\!\!& 
\!(\!\rho^i_j\!\ot\! \sigma\!)[\F\F^{-1}] \!\pm\!
\rho(\Fu\F^{-1(2')}_{(1')})^i_l a_m^+
\sigma(\Fd_{(2)}\F^{-1(2')}_{(2')})
a^l \rho(\F^{(2)}_{(1)}\F^{-1(1')})_j^m
\nonumber \\
 \! &\!\stackrel{(\ref{def4in})}{=}\!&  \! \delta^i_j \idA\pm 
\rho(G^{(1)})_j^m A_m^+\sigma(G^{(2)}) A^l
\rho(G^{(3)})^i_l,
\label{peppe}
\ena 
where
$G:=\F_{12}\F_{23,1}\F^{-1}_{1,23}\F^{-1}_{32}$.
On the other hand, applying to eq. (\ref{cond3})
the permutations $\tau_{23}\circ\tau_{12}$ and
$\tau_{23}$ we obtain respectively
\eqa
\F_{12}\F_{23,1}&=&\F_{31}(\Fu_{(1)} \ot \Fd \ot\Fu_{(2)}) \\
\F^{-1}_{1,23}\F^{-1}_{32}&=&
(\Fu_{(1)} \ot \Fd \ot\Fu_{(2)})^{-1}\F_{13}^{-1},
\ena
which replaced in the definition of $G$ [using
the definition (\ref{rmat})]  give 
$G=\R_{13}$; this proves eq. (\ref{qccr1}).

As for relation (\ref{qccr2}), 
\eqa
A^iA^j &\stackrel{(\ref{def4})}{=}&
\rho(\Fu)^i_l\sigma(\Fd) a^l
a ^m\sigma(\F^{-1(2')})\rho(\gamma^{-1}S \F^{-1(1')})^j_m
\nonumber \\
&\stackrel{(\ref{ccr})}{=}&\pm\rho(\Fu)^i_l\sigma(\Fd) a ^ma^l
\sigma(\F^{-1(2')})\rho(\gamma^{-1}S \F^{-1(1')})^j_m
\nonumber \\
&\stackrel{(\ref{covl})}{=}&\pm\rho(\Fu\F^{-1(2')}_{(2')})^i_l
a ^m\sigma(\Fd_{(1)}\F^{-1(2')}_{(1')}) a^l
\rho(\gamma^{-1}S \F^{-1(1')}\cdot S \Fd_{(2)})^j_m
\nonumber \\
&\stackrel{(\ref{def4in})}{=}&\pm
\rho(\Fu\F^{-1(2')}_{(2')}\F^{-1(\tilde{1})})^i_l
A^m\sigma(\F^{(2'')}\Fd_{(1)}\F^{-1(2')}_{(1')}\F^{-1(\tilde{2})}) 
\nonumber \\
& & \qquad \qquad \qquad \qquad \times
A^l \rho[\gamma^{-1}S (\F^{(1'')}\Fd_{(2)}
\F^{-1(1')}) \gamma]^j_m
\nonumber \\
&\stackrel{(\ref{def1})}{=}&\pm\rho(G^{-1(1)})^i_l 
A^m\sigma(G^{-1(2)}) A^l \rho[S_h(G^{-1(3)}]^j_m.
\ena
But we have already shown that $G=\R_{13}$; by
recalling that
$(\id\ot S_h)\R^{-1}=\R$,  relation (\ref{qccr2})
follows.

Similarly one can prove relation (\ref{qccr3}), which
can be found also more directly by observing 
that in the unitary-\F ~case it follows from the
previous one by applying the $*$-conjugation
and by noting that $\bar R=R_{21}$.
$\Box$

{\it Remark 5.} It is interesting to ask how the 
invariants which can be constructed
from $A^i,A^+_j$ are related to the ones which
can be constructed from $a^i,a^+_j$.
It is straightforward to prove that \eg
any invariant of the form
$I_h^n:=A_{i_1}^+...A^+_{i_n}A^{i_n}...A^{i_1}$
coincides with the invariant of the form
$I^n:=a_{i_1}^+...a^+_{i_n}a^{i_n}...a^{i_1}$
(in particular,
so is $I^2\equiv N:=A_i^+A^i=n $).
One can show \cite{fiocomp} the
equality $P(A^i,A^+_j)=P(a^i,a^+_j)$ for all polynomial
invariants $P(\tilde A^i,\tilde A^+_j)$.
Thi is no more true if $H$ is genuinely quasitriangular.

\sect{Quantum commutation relations: the ${\cal U}_hsl(2)$ case}
\label{sl2cr}

It is now tempting to consider the quasitriangular
case and ask whether a
transformation such as in Remark  1 can map 
$a^i,a^+_j$ into \uqg-covariant
$A^i,A^+_j$ satisfying relations of the type
found in Ref. \cite{puwo,wezu}.
We will stick here to the case $\g=sl(2,\cn)$,
$\rho\equiv$fundamental representation,
addressing the reader to Ref. \cite{fiocomp}
for the general case.

We fix our conventions as follows. 
As `classical' generators of $Usl(2)$
we choose $j_0,j_+,j_-\in\g$ ,
\eqa
[j_0,j_{\pm}] & = &\pm  j_{\pm} \qquad \qquad \qquad \qquad 
  [j_+,j_-] =  2j_0,\\
\Delta (x) & = & {\bf 1}\ot x + x\ot {\bf 1} \qquad \qquad
 \qquad \qquad x=j_0,j_+,j_-;
\ena
as generators of ${\cal U}_hsl(2)$ we choose $J_0,J^+,J^-$ 
satisfying \footnote{Let us recall incidentally that
the mapping $\varphi_h:{\cal U}_hsl(2)\rightarrow Usl(2)[[h]]$,
up to an inner automorphism 
$\varphi_h\rightarrow\varphi_{h,v}:=v\varphi_h(\cdot)v^{-1}$ is given
by \cite{zachos1}
\eq
\varphi_h(J_0)=j_0 \qquad \qquad \qquad 
\varphi_h(J_{\pm})=\sqrt{\frac{[j\pm j_0]_q[1+j\mp j_0]_q}{(j\pm j_0)(1+j\mp j_0)}}
j_{\pm},
\en
where $j$ is the positive root of the equation $j(j+1)-\C=0$,
$\C=j_-j_++j_0(j_0+1)$ is the Casimir, and 
$[x]_q:=\frac{q^x-q^{-x}}{q-q^{-1}}$.}
\eqa
[J_0,J_{\pm}]& = & \pm  J_{\pm} \qquad \qquad \qquad \qquad 
  [J_+,J_-] =  \frac{q^{2J_0}-q^{-2J_0}}{q- q^{-1}}  \label{su2qcr}\\
\Delta_h(J_0) & = & {\bf 1}\ot J_0 + J_0\ot {\bf 1} \qquad \qquad 
\Delta_h(J_{\pm})  =  J_{\pm}\ot q^{-J_0} + q^{J_0}
\ot J_{\pm}.
\label{su2qcocr}
\ena

At the coalgebra level,
the universal \F~ [connecting $\Delta $ to $\Delta_h$ in formula
(\ref{def1})] is not explicitly known; however the
$Usl(2)[[h]]$-valued matrix $F:=(\rho\ot \id)\F$ has been
determined in Ref. \cite{zachos2} and 
reads\footnote{See formulae (3.1), (3.30) in Ref. \cite{zachos2}. To
match our conventions with theirs, one has to rescale 
$j_{\pm}$ by $\sqrt{2}$ and note that the right correspondence
between our notation and theirs
is $\F\equiv\F_q\leftrightarrow U_{q^{-1}}$, what is needed
to match the coproducts.}
\eq
F=\left\Vert
\begin{array}{cc}
        a(j,j_0)  &     b(j,j_0)j_-        \\
-j_+b(j,j_0)  &     a(j,j_0\!-\!1)   
\end{array}
\right\Vert
\qquad\qquad\qquad
F^{-1}=\left\Vert
\begin{array}{cc}
        a(j,j_0)  &     -b(j,j_0)j_-        \\
j_+b(j,j_0)  &     a(j,j_0\!-\!1)   
\end{array}
\right\Vert,
\label{fsu2}
\en
where
\eqa
a(j,j_0) &:=&\frac{q^{j-j_0\over 2}}{\sqrt{(1\!+\!2j)[1\!+\!2j]_q}}
\left[\sqrt{(1\!+\!j\!+\!j_0)[1\!+\!j\!+\!j_0]_q}+ q^{-\frac{(1\!+\!2j)}2}
\sqrt{(j\!-\!j_0)[j\!-\!j_0]_q}\right]  \nonumber \\
b(j,j_0) &:=&\frac{q^{j-j_0\over 2}}{\sqrt{(1+2j)[1+2j]_q}}
\left[\sqrt{\frac{[1+j+j_0]_q}{1+j+j_0}} - q^{-\frac{(1+2j)}2}
\sqrt{\frac{[j-j_0]_q}{j-j_0}}\right];
\label{ffsu2}
\ena
the matrix elements in eq. (\ref{fsu2}) are w.r.t. the
orthonormal basis  $\{|\up\ra,|\downarrow\ra\}$
(spin up, spin down) of eigenstates of $j_0$
with eigenvalues $\pm \frac 12$. 
$a^+_{\up}|0\ra=|\up\ra$ etc. All indices in the sequel
will run over $\{\up,\downarrow\}$. 
$F$ is unitary w.r.t. the $su(2)$ $*$-structure
$j_0^*=j_0$, $(j_+)^*=j_-$.

The homorphism $\sigma$ in this case coincides with the
well-known Jordan-Schwinger realization of $sl(2)$
\cite{bied} and reads
\eq
\sigma(j_+)=a^+_{\up}a ^{\downarrow},\qquad\qquad
\sigma(j_-)=a^+_{\downarrow}a ^{\up},\qquad\qquad
\sigma(j_0)=\frac 12(a^+_{\up}a ^{\up}-
a^+_{\downarrow}a ^{\downarrow})
\label{homo}
\en
implying in particular
\eq
\sigma(j)=\frac {n }2,
\en
where $n :=a_i^+a^i$ is  the `classical 
number of particles' operator.
The ${\cal U}_hsu(N)$-covariant Weyl algebra (with $N\ge 2$,
$q$ real and positive) was introduced by Pusz and Woronowicz
\cite{puwo}; independently, Wess and Zumino \cite{wezu}
introduced its ${\cal U}_hsl(N)$-covariant generalization
(arbitrary complex $q$) in R-matrix notation.
One can consider also its Clifford version \cite{pusz}. 
In the R-matrix notation \cite{wezu}
the QCR of the generators (which we will denote here
by $\tilde A^i,\tilde A^+_j$) read\footnote{The `braiding' (\ref{qccr1n})
between $\tilde A^i,\tilde A^+_j$ could be also replaced by the inverse
one: $\tilde A^i\tilde A^+_j  = \idA \delta^i_j\pm q^{\mp}R^{-1~ui}_{~~~jv}
\tilde A^+_u\tilde A^v$. For a comprehensive introduction to braiding see 
Ref. \cite{majid}}
\eqa
\tilde A^i\tilde A^j          & = & \pm q^{\mp 1}R^{ij}_{vu}\tilde 
A^u\tilde A^v \label{qccr2n}\\
\tilde A^+_i\tilde A^+_j & = & \pm q^{\mp 1}R_{ij}^{vu}
\tilde A^+_u\tilde A^+_v, \label{qccr3n}\\
\tilde A^i\tilde A^+_j     & = &\idA \delta^i_j\pm q^{\pm 1}R^{ui}_{jv}
\tilde A^+_u\tilde A^v ,
\label{qccr1n}
\ena
where $i,j=1,...,N$, $R=(\rho_d\ot \rho_d)\R$
is the $R$-matrix of ${\cal U}_hsl(N)$
in the defining representation $\rho_d$,
and the sign $\pm$ refers to Weyl and Clifford
respectively. Both have been treated subsequently also by many other
authors. They are related, but should not be confused, with the 
celebrated Biedenharn-Macfarlane-Hayashi $q$-oscillator
(super)algebras \cite{mac,hayashi}\footnote{The generators
$\alpha^i,\alpha^+_j$ of the latter fulfil ordinary
(anti)commutation relations, except for the 
$q$-(anti)commutation relations 
$\alpha^i\alpha^+_i \mp q^2\alpha^+_i\alpha^i$ and are {\it not}
$U_h\g$-covariant (in spite of the fact that they are usually
used to construct a generalized Jordan-Schwinger realization
of $U_h\g$). It is of interest to note that, however, the
generators  $\alpha^i,\alpha^+_j$ can be tipically realized
as algebraic `functions' of $\tilde A^i,\tilde A^+_i$
\cite{oleg}, whereas the generators $A^i,A^+_i$ can
be tipically realized only as {\it trascendental} `functions'
of $\tilde A^i,\tilde A^+_i$.}
When $q>0$ the QCR of the generators are compatible with the
$*$-structure \cite{puwo} of annihilators and creators,
\eq
A^+_i=(A^i)^{\star}.
\en
When $N=2$ we will pick up $i,j,...=\up,\downarrow$; the
$R$-matrix will read
\eq
R\equiv\Vert R^{ij}_{hk}\Vert:=\left\Vert
\begin{array}{cccc}
q &                    &      & \\
   &    1              &      &  \\
   &(q-q^{-1}) & 1   &  \\
   &                   &      &q
\end{array}
\right\Vert.
\label{rsu2}
\en
where the row and columns of the matrix (\ref{rsu2}) are ordered
in the usual way: $\up\up, \up \downarrow, \downarrow\up, 
\downarrow\downarrow$ from left to right and from up to down.

More explicitly, the Weyl QCR (\ref{qccr2n}) - (\ref{qccr1n}) 
read
\eqa
\tilde A^{\up}\tilde A^+_{\up}     & = & \idA+q^2 \tilde A^+_{\up}
\tilde A^{\up} +(q^2-1)\tilde A^+_{\downarrow}
\tilde A^{\downarrow}, \nonumber \\
\tilde A^{\downarrow}\tilde A^+_{\downarrow}     
& = & \idA +  q^2 \tilde A^+_{\downarrow}
\tilde A^{\downarrow}, \nonumber \\
\tilde A^{\up}\tilde A^+_{\downarrow}     
& = &  q\tilde A^+_{\downarrow}
\tilde A^{\up}, \nonumber \\
\tilde A^{\downarrow}\tilde A^+_{\up}& = &  q
\tilde A^+_{\up} \tilde A^{\downarrow},
\label{qccr1nnb}
\ena
\eq
\tilde A^{\downarrow}\tilde A^{\up} =  q\tilde A^{\up}
\tilde A^{\downarrow} \qquad\qquad\qquad
\tilde A^+_{\up}\tilde A^+_{\downarrow} =  q
\tilde A^+_{\downarrow}\tilde A^+_{\up} 
\label{qccr2nnb}
\en
and the Clifford
\eqa
\tilde A^{\up}\tilde A^+_{\up}     & = & \idA-\tilde A^+_{\up}
\tilde A^{\up} +(q^{-2}-1)\tilde A^+_{\downarrow}
\tilde A^{\downarrow}, \nonumber \\
\tilde A^{\downarrow}\tilde A^+_{\downarrow}     
& = & \idA - \tilde A^+_{\downarrow}
\tilde A^{\downarrow}, \nonumber \\
\tilde A^{\up}\tilde A^+_{\downarrow}     
& = & - q^{-1}\tilde A^+_{\downarrow}
\tilde A^{\up}, \nonumber \\
\tilde A^{\downarrow}\tilde A^+_{\up}& = & - q^{-1}
\tilde A^+_{\up} \tilde A^{\downarrow},
\label{qccr1nnf}
\ena
\eqa
\tilde A^{\downarrow}\tilde A^{\up} &=&  -q^{-1}\tilde A^{\up}
\tilde A^{\downarrow},\qquad\qquad
\tilde A^{\up} \tilde A^{\up} = 0, \qquad\qquad
\tilde A^{\downarrow}\tilde A^{\downarrow}=0; \\
\tilde A^+_{\up}\tilde A^+_{\downarrow} &=&  - q^{- 1}
\tilde A^+_{\downarrow}\tilde A^+_{\up}, \qquad \qquad
\tilde A^+_{\up} \tilde A^+_{\up} = 0, \qquad\qquad
\tilde A^+_{\downarrow}\tilde A^+_{\downarrow}=0.
\label{qccr2nnf} 
\ena

Now we try to construct the $A^+_i,A^i$.
The Ansatz of Proposition \ref{prop2} is 
equivalent\footnote{In fact, in $T=K[\C\ot \1,\1\ot \C,\Delta (\C)]$
the dependence on the
last argument results only in a numerical factor, because of
eq's (\ref{covl}), (\ref{cov}). The $\C\ot\1$-
and $\1\ot\C$-dependences can be
replaced by the $n$-dependence, since it is easy to
prove that $\sigma(\C)=\frac{n}2(\frac{n}2+1)$; the latter can 
be concentrated
either at the left or at the right of $a^i ,a^+_j$,
upon using the commutation relations 
$[n,a^+_i]=a^+_i$, $[n,a^i]=-a^i$.} to 
\eq
\begin{array}{rcl}
A_i^+ &= &a_l^+\sigma(\F^{-1(2)})\rho(\F^{-1(1)})_i^l f(n )\\
A^i &= &g(n )\rho(\Fu)^i_l\sigma(\Fd) a^l. \\
\label{def4n}
\end{array}
\en
with the same $\F$ and two invertible functions $f,g$.
The product $fg$ is determined by requiring that the commutation
relation between $N:=A^+_iA^i$ and $A^i,A^+_j$ are the
same as those between
$\tilde N:=\tilde A^+_i\tilde A^i$ and $\tilde A^i,\tilde A^+_j$:
\eq
NA^+_i=A^+_i(1+q^{\pm 2}N)\qquad\qquad\qquad
NA^i=A^iq^{\mp 2}(N-1);
\en
one easily finds
\eq
gf(x)=\sqrt{\frac{(x+1)_{q^{\pm 2}}}{(x+1)}}.
\label{funn}
\en
If $q\in\rn^+$ and we wish that $A^+_i=(A^i)^{\star}$,
we need to choose
\eq
g(x)=f(x)=\sqrt{\frac{(x+1)_{q^{\pm 2}}}{(x+1)}}.
\label{funns}
\en
This is nothing else but the already encountered
function (\ref{traproto}) \cite{zachos3} needed to
transform the classical creation/annihilation
operators in one dimension $a^+,a $
into the quantum ones $A^+,A$.

If we are interested in Fock space representations,
eq. $A^+_i|0\ra=a^+_i|0\ra$, stating that the quantum and 
classical one-particle state coincide, is automatically
satisfied, because $f(0)=1$.

Now one can express the RHS(\ref{def4n})
thoroughly in terms of $a^i,a^+_j$. 
It is convenient to introduce the
up-down `number of particle' operators 
$n ^{\up},n ^{\downarrow}$ 
($n ^{\up}+n ^{\downarrow}=n $), by
\eq
n ^{\up}:=a^+_{\up}a^{\up}
\qquad\qquad\qquad
n ^{\downarrow}:=a^+_{\downarrow}
a^{\downarrow}.
\en
Using
Eq.'s (\ref{fsu2}), (\ref{ffsu2}), (\ref{homo}), (\ref{def4n}),
it is easy algebra\footnote{In the fermionic case one has 
to fully exploit the nilpotency of $a^i,a^+_j$.} 
to prove the following 

\begin{prop}
Equations (\ref{def4n}) amount, in the Weyl case, to
\eq
\begin{array}{rclcrcl}
A^+_{\up} & = &\sqrt{(n ^{\up})_{q^2}\over n ^{\up}}
q^{n ^{\downarrow}}a^+_{\up} &\qquad
\qquad A^+_{\downarrow} & = &
\sqrt{(n ^{\downarrow})_{q^2}\over n ^{\downarrow}}
a^+_{\downarrow}  \nonumber \\
A^{\up} & = &a ^{\up}\sqrt{(n ^{\up})_{q^2}\over n ^{\up}}
q^{n ^{\downarrow}} &\qquad
\qquad A^{\downarrow} & = &a ^{\downarrow} 
\sqrt{(n ^{\downarrow})_{q^2}\over n ^{\downarrow}},
\end{array}
\label{lastb}
\en
and in the Clifford one, to 
\eq
\begin{array}{rclcrcl}
A^+_{\up} & = & q^{-n ^{\downarrow}}a^+_{\up} &\qquad
\qquad A^+_{\downarrow} & = &
a^+_{\downarrow}  \nonumber \\
A^{\up} & = &a ^{\up}q^{-n ^{\downarrow}} &\qquad
\qquad A^{\downarrow} & = &a ^{\downarrow}.
\end{array}
\label{lastf}
\en
\end{prop}

We are ready  for the main theorem of this section
(the proof is a straightforward computation).
\begin{theorem}
The elements $A^iA^+_j\in\A$ ($\A$ Weyl or Clifford)
defined in formulae (\ref{def4n}), (\ref{funns}) 
satisfy the QCR (\ref{qccr1n}) - (\ref{qccr3n}).
\end{theorem}

Let us compare the map (\ref{lastb}) (for the Weyl
algebra) with the
one found in Ref. \cite{oleg}. That map, in our
notation, would read
\eq
\begin{array}{rclcrcl}
A^+_{\up} & = & q^{n ^{\downarrow}}a^+_{\up} &\qquad
\qquad A^+_{\downarrow} & = &
a^+_{\downarrow}  \nonumber \\
A^{\up} & = &a ^{\up}{(n ^{\up})_{q^2}\over n ^{\up}}
q^{n ^{\downarrow}} &\qquad
\qquad A^{\downarrow} & = &a ^{\downarrow} 
{(n ^{\downarrow})_{q^2}\over n ^{\downarrow}}
\end{array}
\en
[clearly, it is not compatible with the $*$-structure
$A^+_{i,\alpha}=(A^{i,\alpha})^{\star}$]. It is straightforward to check that
the element $\alpha\in\A[[h]]$ needed to transform $A^+_i,A^i$ into
$A^+_{i,\alpha},A^{i,\alpha}$ [formulae (\ref{fuffi})] is
\eq
\alpha:=\sqrt{\frac{\Gamma(n^{\up}+1)\Gamma(n^{\downarrow}+1)}
{\Gamma_{q^2}(n^{\up}+1)\Gamma_{q^2}(n^{\downarrow}+1)}},
\en
where $\Gamma$ is the Euler $\Gamma$-function and 
$\Gamma_{q^2}$ its $q$-deformation \cite{rahman}
satisfying the property
\eq
\Gamma_q(a+1)=(a)_q\Gamma_q(a).
\en
The corresponding realization of $\tilde\tr$ is obtained through relation
(\ref{goga})

\sect{Representation theory}
\label{repth}

In this section we compare representations of $\A$ with representations
of $\A_h$ and investigate whether the deforming maps found in the
preceding sections can be interpreted as intertwiners between them.

We start with a general remark. At least perturbatively in $h$, we
expect that, for any given representation $\pi$ of $\A$ on some
space $V$, the objects $\pi(A^+_{i,\alpha}),\pi(A^{i,\alpha})$
are well-defined operators on $V$, since
$A^+_{i,\alpha}=a^+_i+O(h)$, $A^{i,\alpha}=a^i+O(h)$;
consequently, $\pi\circ f_{\alpha}$ is a 
representation of $\A_h$ on $V$. Let us prove this statement more
rigorously for some specific kind of representations.

If $\g$ is compact, we can choose $\rho$ to be a unitary representation;
the $\g$-covariant CCR will admit the $*$-relations (\ref{star1}).
Assuming the latter, Stone-Von Neumann theorem (or its Clifford
counterpart) applies: there exists a `$*$-representation' $\pi$
of eq. (\ref{ccr}) and (\ref{star1}) on a separable Hilbert space
$\Hil$, \ie a $\pi$ fulfilling the following properties:
\begin{enumerate}
\item $\pi(a^i)$, $\pi(a^+_i)$ are closed;
\item $\pi(a^+_i)\subset[\pi(a^i)]^{\dagger}$;
\item there exists a dense linear subset ${\cal D}\subset\Hil$
contained in the domain of the product of any two operators 
$\pi(a^i)$, $\pi(a^+_i)$;
\item the CCR (\ref{ccr}) hold on  ${\cal D}$;
\item the nonnegative-definite operator 
$\pi(n):=\pi(a^+_ia^i)$ is essentially self-adjoint on ${\cal D}$.
\end{enumerate}

Moreover, the irreducible components $(\pi_l,\Hil_l)$ of $(\pi,\Hil)$ are,
up to a unitary transformation $U_l:\Hil\rightarrow \Hil$,
Fock space representations, \ie there exist `ground states'
$|0\ra_l\in\Hil_l$:
\eq
\pi(a^i)|0\ra_l={\bf 0}.
\en
One can choose the dense set $\cal D$ as the linear span of all
the analytic vectors of the form 
$\pi(a^+_{i_1}...a^+_{i_k})|0\ra_l$.

In the case under consideration, choosing $\F$ unitary [so that relation
(\ref{star2}) hold] and setting $\Pi:=\pi\circ f$, one can easily realize
that $\Pi$ is a $*$-representation of the QCR [(\ref{QCR0}) or 
(\ref{qccr2n}-\ref{qccr1n})] on
$\Hil$, \ie it fulfils conditions analogous to 1. - 5. ,
with $\tilde A^i,\tilde A^+_i$ in the place of $a^i,a^+_i$.

The reason is essentially the following. Assume first that $\pi$ and $\rho$
are irreducible, and let 
$\Hil^{(k)}:=\mbox{Span}_{\cn}\{\pi(a^+_{i_1}...a^+_{i_k})|0\ra\}$;
$\Hil^{(k)}$ is a finite dimensional eigenspace of $\pi(n)$, because
$\rho$ is finite-dimensional and therefore, $\A$ has a finite
number of generators. From equations (\ref{def4}) and the relation
$[\sigma(U\g),n]=0$ it follows that $\pi(A^i),\pi(A^+_i)$  differ from
$\pi(a^i),\pi(a^+_i)$ just by operators mapping each $\Hil^{(k)}$
into itself\footnote{Incidentally, from relation (\ref{cond2}) it 
follows in particular $\pi(A^i)|0\ra={\bf 0}$ and 
$\pi(A^+_i)|0\ra=\pi(a^+_i)|0\ra$.}. Therefore all the properties
1. - 5. are inherited by $\pi(A^i),\pi(A^+_i)$ as well.
The same result holds for $\Pi_{\alpha}:=\pi\circ f_{\alpha}$
if $\alpha^{\star}=\alpha^{-1}$, because $\pi(\alpha)$
can be absorbed in the unitary transformation $U$. Finally, 
the result extends by linearity and orthogonality also to the case
that $\pi$ and/or $\rho$ are reducible.

Summing up, $f_{\alpha}$ is an intertwiner between the 
$*$-representation of the CCR and a $*$-representation of the QCR.
In other words, we can represent the objects 
$\tilde A^i,\tilde A^+_i$ [fulfilling the QCR and the $*$-relations
(\ref{star2})] as composite operators acting on the Fock space
$\Hil$, as anticipated in the introduction. The latter can be used
to describe ordinary Bosons and Fermions.
This disproves the quite common belief that non-cocommutative Hopf
algebra symmetries are necessarily incompatible with ordinary
Bose and Fermi statisitcs.

Let us analyze now a different situation. Let $\g=sl(2)$,
$\A$ be the corresponding two-dimensional Weyl algebra
and $(\pi,V)$ the Fock space representation of the latter
(with ground state $|0\ra$); let $\A_h$ be the deformation of $\A$
considered in section \ref{sl2cr} with $q^2=e^{2h}$ a root of unity,
\ie $q^{2p}=1$, $q^{2k}\neq 1$ with $p,k\in\nn$ and $k<p$. It is
easy to realize that $\pi[(n_i)_{q^2}(a^+_i)^{mp}(a^+_{-i})^l]|0\ra={\bf 0}$
for $l,m=0,1,...$; consequently, $A^{\up},A^{\downarrow}$ annihilate
all the vectors of the form 
$|m,n\ra:=\pi[(a^+_{\up})^{mp}(a^+_{\downarrow})^{np}|0\ra$
($n,m=0,1,...$). Although $(\pi,V)$ was irreducible as a
representation of $\A$, it is reducile as a representation
of the subalgebra of $\A$ generated  by $A^i,A^+_j$; the irreducible
components $V_{m,n}$ are isomorphic  and $p^2$-dimensional,
and are obtained by applying $A^+_i$'s to the cyclic vectors
$|m,n\ra$. $V_{m,n}$ is also a (reducible) representation of 
${\cal U}_hsl(2)$ and may be called (with an abuse of terminology,
since we have not introduced any spatial degrees of freedom)
an `anyonic space'. Thus, $f$ can be seen as an intertwiner
from the classical Bosonic Fock representation onto a direct
sum of `anyonic' representations of $\A_h$.

Whenever some class $\sf P$ of representations of 
$\A_h$ is `larger' than the corresponding class $\sf p$ of
representations of the CCR, then $f_{\alpha}^{-1}$ can of
course be well-defined (as an intertwiner between $\sf P$ and
$\sf p$) only on some proper subset of $\sf P$; on its
complement it must be singular. 

It is instructive to see whether and how this phenomenon occurs in some
concrete example, \eg for the class of $*$-representations considered
above. Both a deformed algebra $\A_h$ covariant under some triangular
Hopf algebra $H_h$ of the type (\ref{reshet}) (with $\F$ unitary) and the
deformed Clifford algebra of section \ref{sl2cr} are not good for this
purpose, because $\sf P$ is as large as 
$\sf p$\footnote{For the latter algebra this was shown in Ref.
\cite{pusz}. For the former this can be understood as follows.
In the present case $\R=\F^{-2}$. Given any representation $\rho$
of $U\g$ and the corresponding representation $\tilde \rho$
of $U_h\g$, let us choose a basis of eigenvectors $|l\ra$
of the Cartan subalgebra. In this basis 
$R:=(\tilde \rho\ot \tilde \rho)\R$ will be diagonal, and in
particular it will be $R^{ll}_{ll}=1$ as a
consequence of the antisymmetry of $\omega_{ij}(h_i\ot h_j)$.
The QCR (\ref{QCR0}) will imply in particular
\eq 
\tilde A^i\tilde A^+_i=\idA\pm \tilde A^+_i\tilde A^i
\z\z\z (\mbox{no sum over }i),
\en
and, setting $\tilde N_i:=A^+_i\tilde A^i$
(no sum over $i$),
\eq
[\tilde N_i,\tilde N_j]=0\z\z
[\tilde N_i,\tilde A^+_j]=\delta_{ij}\tilde A^+_j,
[\tilde N_i,\tilde A^j]=-\delta_i^j\tilde A^j.
\label{prora}
\en
Let us denote by $\A_h{}^{(i)}$ the subalgebra generated by 
$\tilde A^i,\tilde A^+_i$. Each $\A_h{}^{(i)}$ separatly is isomorphic
to a classical one-dimensional Weyl/Clifford algebra, and therefore
admits, up to unitary equivalences, a unique $*$-representation in the
form of a Fock representation with level $N_i$. Because of the
relations (\ref{prora}), we can choose $\{N_i\}_{i\in I}$ as a
complete set of commuting observables, whence the uniqueness, up to
unitary equivalences, of the $*$-representation of $\A_h$ follows
immediately.}. 
On the contrary,
it was proved in Ref. \cite{puwo} that there are many unitarily
non-equivalent $*$-representations on separable Hilbert spaces of 
the ${\cal U}_hsu(2)$-covariant ($q\in\rn^+$) deformed
Weyl algebra\footnote{Incidentally, even 1-dimensional deformed
Heisenberg algebras may have more unitarily inequivalent representations
\cite{wess}. Moreover, within each representation one has
still some freedom in the `physical' interpretation of the
observables, \eg what are the `right' momentum/position observables,
see \eg Ref. \cite{kempf}.}.

Sticking to the case $0<q<1$, one can parametrize
the Woronowicz-Pusz \cite{puwo} unitarily
inequivalent irreducible $*$-representations
of this algebra by three parameters $E,r,s$, where $q^2<E<1$, 
and $r,s$ are nonnegative integers with $r+s\le 2$. 
We shall denote the corresponding Hilbert space by
$\Hil_{E,r,s}$.
We divide them in the following classes for clarity:

\begin{enumerate}

\item In the representation 
$s=2$ (and $r=0$, the value $E$ is irrelevant) one 
parametrizes the vectors of an orthonormal basis of $\Hil_{E,0,2}$
 by
$\{|\frac{q^{2m_1}}{q^2-1},\frac{q^{2m_2}}{q^2-1}\ra\}$.

\item In the representations with $s=1=r$, one 
parametrizes the vectors of an orthonormal basis of $\Hil_{E,1,1}$
 by
$\{|q^{2n_1}E,\frac{q^{2m_1}}{q^2-1}\ra\}$.

\item In the representations with $s=1$, $r=0$, one 
parametrizes the vectors of an orthonormal basis of $\Hil_{E,0,1}$
by
$\{|0,\frac{q^{2m_1}}{q^2-1}\ra\}$
(the value $E$ is irrelevant).

\item In the representations with $s=0$, $r=2$, one 
parametrizes the vectors of an orthonormal basis of $\Hil_{E,2,0}$
by
$\{|q^{2n_1}E,q^{2n_2}E\ra\}$.

\item Finally, in the representations with $s=0$, $r=1$, one 
parametrizes the vectors of an orthonormal basis of $\Hil_{E,1,0}$ by
$\{|q^{2n_1}E,0\ra\}$.

\end{enumerate}

Here $n_1,n_2;m_1,m_2$ denote integers, with 
$m_1\ge m_2\ge 0$, $n_1<n_2$. Only representations 1, 3 
have a ground state; but representation 3 is 
degenerate. 

In the rest of this section we drop the symbols $\Pi$ to
avoid a too heavy notation.
On the vectors 
$|\eta_{\up},\eta_{\downarrow}\ra$ of the above
basis of any $\Hil_{E,r,s}$,
$\tilde A^i,\tilde A_i^+$, $i=\up,\downarrow$,
are defined (modulo a possible but here irrelevant
phase in the case $r+s<2$) by  
\eq
\begin{array}{rclcrcl}
\tilde A^{\up}|\eta,\eta_{\downarrow}\ra &=&
\sqrt{\eta -\eta_{\downarrow}}|q^{-2}\eta ,
\eta_{\downarrow}\ra
\x\qquad &
\tilde A^{\downarrow}|\eta ,\eta_{\downarrow}\ra &=&
\sqrt{\eta_{\downarrow} -{1\over q^2-1}}|q^{-2}
\eta ,q^{-2}\eta_{\downarrow}\ra \\
\tilde A^+_{\up}|\eta ,\eta_{\downarrow}\ra &=&
\sqrt{q^2\eta -\eta_{\downarrow}}|q^2\eta ,\eta_{\downarrow}\ra
\x\qquad &
\tilde A^+_{\downarrow}|\eta ,\eta_{\downarrow}\ra &=&
\sqrt{q^2\eta_{\downarrow}-{1\over q^2-1}}|q^2\eta ,
q^2\eta_{\downarrow}\ra. \\
\end{array}
\label{genrep}
\en
Hence, setting $\tilde N^{\up}:=\tilde A^+_{\up} \tilde A^{\up}$,
$ \tilde N^{\downarrow}:=\tilde A^+_{\downarrow} \tilde A^{\downarrow}$
and $\tilde N:=\tilde N^{\up}+\tilde N^{\downarrow}$,
 we find that $|\eta,\eta_{\downarrow}\ra$ are 
eigenvectors of the following operators:
\eq
\begin{array}{lll}
\tilde N^{\up}|\eta ,\eta_{\downarrow}\ra &=&
(\eta -\eta_{\downarrow})|\eta ,\eta_{\downarrow}\ra\cr
\tilde N^{\downarrow}|\eta ,\eta_{\downarrow}\ra & =&
(\eta_{\downarrow} -{1\over q^2-1})|\eta ,
\eta_{\downarrow}\ra,\cr
[1+(q^2-1)\tilde N^{\downarrow}]|\eta ,\eta_{\downarrow}\ra &=&
(q^2-1)\eta_{\downarrow}|\eta,\eta_{\downarrow}\ra,\cr
[1+(q^2-1)\tilde N]|\eta ,\eta_{\downarrow}\ra &= &
(q^2-1)\eta|\eta,\eta_{\downarrow}\ra.
\end{array}
\label{corre}
\en

Formally, the inverse of
the transformation (\ref{lastb}) reads
\eqa
f^{-1}(\tilde A^+_{\downarrow}):=\tilde a^+_{\downarrow}&=&\sqrt{\log[1+(q^2-1)\tilde 
N^{\downarrow}]
\over 2\tilde N^{\downarrow} \log q}\tilde A^+_{\downarrow}
\nonumber\\
f^{-1}(\tilde A^+_{\up}):=\tilde a^+_{\up}&=&\sqrt{\frac{1+(q^2-1)\tilde
 N^{\downarrow}}
{2\tilde N^{\up}\log q}\log\left[\frac{1+(q^2-1)\tilde N}
{1+(q^2-1)\tilde N^{\downarrow}}\right]}\tilde A^+_{\up}
\nonumber\\
f^{-1}(\tilde A^{\downarrow}):=\tilde a ^{\downarrow}&=&\tilde A^{\downarrow}
\sqrt{\log[1+(q^2-1)\tilde 
N^{\downarrow}]\over 2\tilde N^{\downarrow} \log q}
\label{invtran} \\
f^{-1}(\tilde A^{\up}):=\tilde a ^{\up}&=&\tilde 
A^{\up}\sqrt{\frac{1+(q^2-1)\tilde N^{\downarrow}}
{2\tilde N^{\up}\log q}\log\left[\frac{1+(q^2-1)\tilde N}
{1+(q^2-1)\tilde N^{\downarrow}}\right]}.
\nonumber
\ena
A glance at formula (\ref{corre}) shows that the arguments
of both logarithms
in the inverse transformation (\ref{invtran}) are
positive-definite on representation 1, whereas at least one of the
arguments of the two logarithms is negative on any
representaion 2,3,4 or 5 (in the representation 3 the
other argument vanishes). This makes $f^{-1}$ ill-defined
on all representations 2,3,4 or 5.

We conclude that  representation 1 is the one intertwined
by $f^{-1}$ to the 
standard bosonic Fock representation of the
$su(2)$-covariant Weyl algebra $\A$, whereas
the representations of the classes
2,4,5 have no classical analog, and representation 3
reduces to the representation of a 1-dimensional
Weyl algebra.

\section*{Acknowledgments}

I would like to thank P. Schupp for fruitful discussions on the subject. 
I am grateful to J.~Wess for the 
warm hospitality at his institute. This work was financially supported
by a TMR fellowship grranted by
the European Commission, Dir. Gen. XII for Science,
Research and Development, under the contract ERBFMBICT960921.

\end{document}